\documentstyle[11pt]{article}
\newtheorem {th}{Theorem}[section]
\newtheorem {lem}[th]{Lemma}
\newtheorem {pr}[th]{Proposition}
\newtheorem {cor}[th]{Corollary}

\def\Cox{\hfill \Box}

\def\ul{\underline}

\def\ee{\epsilon}

\def\E{{\bf{E}}}

\def\P{{\bf{P}}}

\def\Z{{\bf Z}}
\def\F{{\cal{F}}}
\def\|{\, | \, }
\def\one{{\bf 1}}

\def\abcd{{a\;\;b \choose c \;\; d}}

\def\Var{{\rm Var}}

\begin{document}

\begin{titlepage}
\begin{center}
{\large \bf Vertex-reinforced random walk on $\Z$ has finite range} \\
\end{center}
\vspace{5ex}
\begin{flushright}
Robin Pemantle \footnote{Research supported in part by a Sloan Foundation
Fellowship, and by a Presidential Faculty Fellowship}$^,$\footnote{Department 
of Mathematics, University of Wisconsin-Madison, Van Vleck Hall, 480 Lincoln
Drive, Madison, WI 53706}
  ~\\
Stanislav Volkov \footnote{
The Fields Institute for Research in Mathematical Sciences,
222 College Street, Toronto, Ontario, Canada, M5T3J1}
\end{flushright}

\vfill

{\bf ABSTRACT:} \break
A stochastic process called Vertex-Reinforced Random Walk (VRRW) is 
defined in Pemantle (1988a).  We consider this process in the case
where the underlying graph is an infinite chain (i.e., the one-dimensional
integer lattice).  We show that the range is almost surely finite,
that at least 5 points are visited infinitely often almost surely,
and that with positive probability the range contains exactly 5 points.
There are always points visited infinitely often but at a set of times
of zero density, and we show that the number of visits to such a point
to time $n$ may be asymptotically $n^\alpha$ for a dense set of 
values $\alpha \in (0,1)$.  The power law analysis relies on analysis
of a related urn model.
\vfill


\end{titlepage}

\setcounter{equation}{0}

\section{Outline of results}

For any process $X_0 , X_1 , X_2 , \ldots$ taking values in the vertex
set of a graph $G$ (throughout this paper $G=\Z$), we define the augmented occupation numbers 
$$Z(n,v) = 1 + \sum_{i = 0}^n \one_{X_i = v}$$
to be the number of times plus one that the process visits site $v$
up through time $n$.  Let $G$ be any locally finite graph, with
the neighbor relation denoted by $\sim$, and define vertex-reinforced 
random walk (VRRW) on $G$ with starting point $v \in V(G)$ 
to be the process $\{ X_i : i \geq 0 \}$ such that $X_0 = v$ and
$$\P (X_{n+1} = x \| \F_n) = \one_{x \sim X_n} {Z(n,x) 
   \over \sum_{w \sim X_n} Z(n,w)} \, .$$
In other words, moves are restricted to the edges of $G$, with the
probability of a move to a neighbor $w$ being proportional to the 
augmented occupation of $w$ at that time.  

This is a special case of the weighted VRRW, defined by 
Pemantle (1988a, 1992), where each oriented edge $\vec{vw}$ carries a
nonnegative weight $\lambda (v,w)$, and the transition
probabilities are given by
$$\P (X_{n+1} = x \| \F_n) = \one_{x \sim X_n} {\lambda (X_n , x) Z(n,x) 
   \over \sum_{w \sim X_n} \lambda (X_n , w) Z(n,w)} \, .$$
It is shown in Pemantle (1988a, 1992) that for generic 
symmetric
values of $\lambda$ and finite graphs $G$, the vector of normalized occupation
measure, $(Z(n,v)/n))_{v \in V(G)}$, must converge to an element of
a set of equilibrium points which is typically finite.  Unfortunately,
the case where $\lambda$ is identically 1 is not generic, but rather
degenerate from the point of view of the previous works, and so this
one most natural case is left unanalyzed.  

While the results of Pemantle (1992) do not extend to the case in this
paper, it was conjectured there that in such cases the range of
VRRW will be finite, and that in fact it will get ``stuck'' in 
a set of 3 points.  In this paper we show that this behavior holds,
in the sense of normalized occupation measure, at least with 
positive probability (Theorem~\ref{th 3} below).  If one cares
about the set of points visited infinitely often, rather than with 
positive density, the size of the set on which the walk gets 
stuck is 5 rather than 3.  We obtain the following further results 
about the size of the range.

\begin{th} \label{th 1}
Let $R = \{ k : X_n = k \mbox{ for some $n$} \}$ be the (random) range of 
the process $X_0 , X_1 , \ldots$.
Then $\P (|R| = 5) > 0$ and $\P (|R| < \infty) = 1$.
\end{th}

\begin{th} \label{th 2}
Let $R' = \{ k : X_n = k \mbox{ infinitely often} \}$ be the
essential range of the process $X_0 , X_1 , \ldots$.  Then 
$\P (|R'| \leq 4) = 0$.
\end{th}

\noindent{\em Remark:} Simulations appear to show that $\P (|R'| = 4)$
is nonzero, but these are evidently misleading.

We conjecture but cannot prove that $\P (|R'| = 5) = 1$.
It is easy to see that if $R' = \{ k , k+1 , \ldots , k+j \}$,
then $Z(n,k) / n$ and $Z(n,k+j) / n$ both converge to zero, or in 
other words, the occupation density goes to zero at the endpoints
of the range.  Quantitatively, we have

\begin{th} \label{th 3}
For any closed interval $I \subseteq (0,1)$ and any integer $k$,
there is with positive probability an $\alpha \in I$ such that
the following six events occur: 
\begin{quote}
$(i)$ $R' = \{ k-2 , k-1 , k , k+1 , k+2 \}$; \\
$(ii)$ $\log Z(n , k+2) / \log n \rightarrow \alpha$; \\
$(iii)$ $\log Z(n , k-2) / \log n \rightarrow 1 - \alpha$; \\
$(iv)$ $Z(n , k+1) / n \rightarrow \alpha / 2$; \\
$(v)$ $Z(n , k-1) / n \rightarrow (1 - \alpha) / 2$; \\
$(vi)$ $Z(n , k) / n  \rightarrow 1/2$.
\end{quote}
\end{th}
We conjecture but cannot prove that this is the universal behavior, 
i.e., that there is always such an $\alpha \in (0,1)$.  

The power law behavior in parts~$(ii)$ and~$(iii)$ of Theorem~\ref{th 3} 
rests on the analysis of a certain interacting urn process.  This urn
process is of a type studied in the doctoral dissertation of Athreya
(1967), via embedding in a multiptype branching process.  Since this
is not generally available, and since our hypotheses and methods of proof
are quite different, we include complete statements and proofs of the
relevant results.  The next section gives some background on urn processes
and reinforced random walks.  Proofs for the results on urns are given in 
Section~3, and proofs of Theorems~\ref{th 1},~\ref{th 2}, 
and~\ref{th 3} are given in Section~4.  The final section completes
the proofs of some lemmas and poses a few open questions.

\setcounter{equation}{0}

\section{Background on urn processes and processes with reinforcement}

This section begins with a brief survey of previously studied reinforced 
random processes.  By popular demand, we have included more of a review
than is strictly necessary for the analysis of VRRW, that being the
generalization of Theorem~\ref{th urn 1} stated and proved in the next
section.

The simplest (and one of the oldest) process with reinforcement is
known as {\em P\'olya's urn}, after the 1927 paper of Eggenberger
and  P\'olya.  In this model, there is an urn containing red and 
blue balls.  At time 0 the urn contains $r$ red balls and $s$ blue balls.
At each time $k \geq 1$, a ball is chosen uniformly from the contents 
of the urn, and is put back into the urn along with $a$ extra balls
of the same color.  Thus if $X_n$ denotes the number of red balls
at time $n$ and $Y_n$ denotes the number of blue balls at time $n$,
the dynamics are governed by
\begin{eqnarray} \label{eq polya} 
(X_{n+1} , Y_{n+1}) = (X_n + a , Y_n) && \mbox{ with probability } 
   {X_n \over X_n + Y_n} ; \\[1ex]
(X_{n+1} , Y_{n+1}) = (X_n , Y_n + a) && \mbox{ with probability } 
   {Y_n \over X_n + Y_n} . \nonumber 
\end{eqnarray}
Eggenberger and P\'olya showed that the proportion of red balls,
$Z_n := X_n / (X_n + Y_n)$, converges almost surely, and that the limit
is random.  The distribution of the limit is a beta with parameters 
$r/a$ and $s/a$ (thus uniform over $[0,1]$ when $r=s=a=1$).  This
random limit behavior is possible because $Z_n$ is 
a martingale.  In the sections to follow, we make use several times 
of this elementary principle:
\begin{pr} \label{pr polya}
If $Z_n = X_n / (X_n + Y_n)$, and if $(X_{n+1} , Y_{n+1}) = (X_n + 1 , Y_n)$
with probability $X_n / (X_n + Y_n)$, and $(X_n , Y_n + 1)$ otherwise,
then $\E (Z_{n+1} \| Z_n) = Z_n$.  Also, $|Z_{n+1} - Z_n| < 1 / (X_n + Y_n)$.
$\Cox$
\end{pr}

In contrast to this is an variant suggested by B. Friedman (1949),
where in addition to the $a$ extra balls of the same color, one
also adds $b$ balls of the opposite color.  This produces strikingly
different behavior, even when $b << a$.  Freedman (1965) showed that
$Z_n \rightarrow 1/2$ almost surely, with $(Z_n - 1/2) / n^{- \gamma}$
converging to a nontrivial distribution for $\gamma > 0$ depending
on $a$ and $b$.  To explain the differing behavior, note that
$\{ Z_n \}$ is not a martingale, but rather 
\begin{equation} \label{eq friedman}
\E (Z_{n+1} - Z_n \| Z_n) = n^{-1} (f(Z_n) + o(1))
\end{equation}
where $f$ is a function vanishing only at $1/2$.  One could say that
the drift, $f$, pushes $\{ Z_n \}$ toward $1/2$, which is an 
attracting point for the one-dimensional vector field given by $f$.  

When discussing processes with reinforcement, it is good to keep in
mind the distinction {\em P\'olya-like} ($f \equiv 0$) 
versus {\em Friedman-like} ($f \neq 0$ except at isolated points),
which dictates the important aspects of the long-term behavior.
A third category, {\em singular}, occurs when $f$ has zeros on
the boundary, in which case convergence happens at a slower rate.

The prototypical Friedman-like model is Robbins and Monro's (1951)
stochastic approximation scheme, which obeys the law 
$$\E (Z_{n+1} - Z_n \| Z_n) = n^{-1} F(Z_n)$$
for a generic function $F$ about which imprecise information can
be obtained by sampling.  The (unknown) zeros of $F$ are then
``found'' by the $\{ Z_n \}$ process.  Since the 1950's, stochastic
approximation has been an active research area; the overview by Kushner 
and Yin (1997) gives an idea of progress and techniques in stochastic 
approximation since then.  The literature on formal
models of learning contains many Friedman-like processes, in which
the transition probabilities of a finite state (non-Markov) chain 
are updated based on the some kind of objective function.  The 
chain then ``learns'' to spend most of its time at states with
large values of the objective function.  The first round of this
literature appeared in the late sixties, e.g., Iosifescu and Theodorescu
(1969), and a second round emerged with the study of neural nets.
The common theme is self-organization by a system whose basic
parameters are extremely simple.  

P\'olya-like models have appeared frequently in theoretical statistics,
due to the fact that bounded martingales are mathematically equivalent 
to sequences of posteriors, given increasing $\sigma$-fields.  For example,
suppose an IID sequence of zeros and ones has an unknown mean $p$, with
the prior on $p$ being uniform on $[0,1]$ (or more generally, any beta
distribution).  Then the sample sequences $\{ Z_n \}$ of a 
P\'olya urn process can be interpreted as posterior means, where
each red ball picked corresponds to observing a one and each blue ball
picked corresponds to observing a zero.  The so-called Bayes-Laplace
estimate of the probability the sun will rise tomorrow and 
Greenwood and Yule's (1920) model for industrial accidents are
both based on this interpretation.  In modern times, Blackwell
and McQueen (1973) contruct Ferguson's Dirichlet via an urn process,
and Mauldin, Sudderth and Williams (1992) use a tree full of urns to construct
a family of priors on distributions on $[0,1]$ with nice properties.

P\'olya-like models have also arisen in modeling of self-organization
and random limits.  Arthur (1986) and Arthur et al (1987) 
use both P\'olya- and Friedman-like urns to model the growth of 
industry and explain random clustering and market share patterns.
Reinforced random walks were introduced by Coppersmith and
Diaconis (1987) as another, somewhat simplified model of 
self-organized behavior.  Although simplistic, urn models and
reinforced random walks have been taken seriously in the modeling
of physical phenomena; see for example Othmer and Stevens (1998),
in which motion and aggregation of myxobacteria along slime trails
are modeled by reinforced random walks and related stochastic cellular
automata.  

The VRRW studied by Pemantle (1988a , 1992)
is a variant of their edge-reinforced random walk (ERRW).  In ERRW, one
keeps track of the number of times each edge has been crossed,
the augmented occupation numbers being denoted $\{ Z(n , \{ v , w \}) 
: \{ v , w \} \in E(G) \}$, and one chooses the next edge from among 
the edges adjacent to the present vertex, with probabilities 
proportional to the augmented occupation of each edge:  
$$\P (X_{n+1} = x \| \F_n) = \one_{x \sim X_n} {Z(n , \{ X_n , x \}) 
   \over \sum_{w \sim X_n} Z(n , \{ X_n , w \} ) } \, .$$
Reinforcing edges rather than vertices makes a dramatic difference in
the behavior of the process, because edge-reinforcement is 
P\'olya-like and vertex-reinforcement is Friedman-like.  
A depiction of this difference via simulation may be found in Othmer 
and Stevens (1998).  To see how to account for the difference 
theoretically, let $v$ is a vertex in an acyclic graph, with
incident edges $e_1 , \ldots , e_k$.  The successive edges chosen
each time $v$ is visited form a P\'olya urn process and it is not hard 
to see that these are independent as $v$ ranges over all vertices.
(The analogue of this fact on a graph with cycles is much harder to
formulate and prove.)
Coppersmith and Diaconis (1987) proved that the normalized occupation
measure of ERRW 
on a finite graph 
converges to a random vector having a nonzero density with
respect to Lebesgue measure on the simplex.  When $G = \Z$, 
Pemantle (1988b) shows that the process is a mixture of positive 
recurrent Markov chains, and in particular, that the normalized 
occupation measure converges to a limit that is everywhere positive.  

The question of the behavior of either VRRW or ERRW on a lattice
of dimension two or greater is still open.  Some progress on ERRW
has been made by generalizing the model so that the $k^{th}$ crossing
of each edge adds $a_k$ to the occupation, where $\{ a_k \}$ is
a pre-specified sequence ($a_k \equiv 1$ in standard ERRW).  A
general recurrence/transience dichotomy for this model was
obtained by Davis (1990) in one dimension, while Sellke (1994)
has results on the coordinate processes for this model in two dimensions.

Our results for one-dimensional VRRW depend on an analysis of an
urn model generalizing both the P\'olya and the Friedman urn.  
Replace the dynamics~(\ref{eq polya}) by the more general dynamics:
\begin{eqnarray} \label{eq urn}
(X_{n+1} , Y_{n+1}) = (X_n + a , Y_n + b) && \mbox{ with probability } 
   {X_n \over X_n + Y_n} ; \\[1ex]
(X_{n+1} , Y_{n+1}) = (X_n + c , Y_n + d) && \mbox{ with probability } 
   {Y_n \over X_n + Y_n} . \nonumber 
\end{eqnarray}
There is no assumption that the number of balls be integral.  When
$\abcd$ is a multiple of the identity matrix, we recover P\'olya's
urn, and when $a=d$ and $b=c$ are all nonzero, we recover Friedman's urn.  
In any case where $\abcd$ has an eigenvector $(v_1 , v_2)$ with 
positive components, Freedman's analysis can be carried through 
to show that $X_n / (X_n + Y_n)$ converges to $v_1 / (v_1 + v_2)$.
Perhaps the cleanest way to do this is via embedding in a branching
process, as described in Athreya and Ney (1972, chapter V, sec. 9).
Thus in particular this holds when $bc > 0$.  Two interesting cases
are the singular cases, which can be reduced without loss of generality
to the cases in the next two theorems.  Theorem~\ref{th urn 1}
was first proved by Athreya (1967) in a different form, while
Theorem~\ref{th urn 2} is derivable from his results.

\begin{th} \label{th urn 1}
Suppose $a > d = 1$ and $b = c = 0$.  Then $X_n / Y_n^a$ converges
almost surely to a random limit in $(0,\infty)$.   
\end{th}

\begin{th} \label{th urn 2}
Suppose $a = d = 1, b = 0$ and $c > 0$.  Then $X_n / (c Y_n) - \log Y_n$
converges to a random limit in $(-\infty , \infty)$.
\end{th}

\noindent{\em Remarks:} (1) Theorem~\ref{th urn 2} is in a sense a 
finer result than Theorem~\ref{th urn 1}, since it deals with the
second order correction: $Y_n$ is like $n / \log n$ multiplied
by a specific constant, with a random correction of lower order:
$X_n \approx c Y_n (A + \log Y_n)$, where $A$ is random. 
(2) The class of urns in Theorem~\ref{th urn 2} is not needed 
for analysis of VRRW on $\Z$, but is relevant to VRRW for a different
reason.  In the case $c = 1$, there is an isomorphism between the 
urn process and VRRW on the graph $G$ with $V(G) = \{ A , B \}$, 
having one edge between $A$ and $B$ and one loop connecting $A$ to itself.  
Thus VRRW on $G$ spends roughly time $n / \log n$ at $B$ up to time $n$.

\setcounter{equation}{0}

\section{Urn model proofs}

This section is devoted to proving Lemma~\ref{lem general}, which
generalizes Theorem~\ref{th urn 1} to allow random increments. 
Whereas Athreya (1967) proved version of this by embedding in a multitype
branching process, we use martingale arguments (also considered by
Athreya in some subcases).  These turn out to
be easier in the case of Theorem~\ref{th urn 2} than in the case of 
Lemma~\ref{lem general} below.  Consequently, we first give a relatively 
short proof of Theorem~\ref{th urn 2} and then state and prove
Lemma~\ref{lem general}.  Depending on your tastes, you may find 
the shorter proof or the more modular general proof easier to follow.
Begin with the following easy lemma.
\begin{lem} \label{lem to infty}
Let the nonnegative matrix $\abcd$ satisfy $(a + c) (b + d) > 0$ 
and define an urn process as in~(\ref{eq urn}).  Then $\min 
\{ X_n , Y_n \} \rightarrow \infty$ almost surely.
\end{lem}

\noindent{\sc Proof:} The proportion of red balls at time $n$
is always at least $X_0 / (X_0 + Y_0 + n(a+b+c+d))$.  Since the
sum of these quantities is infinite, the Borel-Cantelli Lemma tells
us that a red ball is chosen infinitely often.  Similarly, a blue
ball is chosen infinitely often.  After each color has been
chosen $k$ times, $\min \{ X_n , Y_n \}$ is at least
$k \min \{ a+c , b+d \}$.   $\Cox$

The following general fact about convergence of random sequences is
also useful.
\begin{lem} \label{lem semi}
Let $\{ Z_n : n \geq 0 \}$ be a random sequence measurable with respect
to the filtration $\{ \F_n \}$.  Define
$$\Delta_n = \E (Z_{n+1} - Z_n \| \F_n ) \;\; \mbox{ ; } \;\;
   Q_n = \E ((Z_{n+1} - Z_n)^2 \| \F_n ) .$$
Then as $n$ goes to infinity, $Z_n$ converges to a finite value almost
surely on the event $\sum_n \Delta_n < \infty$ and $\sum_n Q_n < \infty$.
\end{lem}

\noindent{\sc Proof:} Let $\tau_M$ be the first time $n$ that 
$\sum_{j=0}^n Q_j > M$.  Let 
$$Z_n^{(M)} = Z_{n \wedge \tau_M} - \sum_{j=0}^{n \wedge \tau_M} 
   \Delta_j \, .$$
Observe that $\{ Z_n^{(M)} \}$ is a martingale with 
$$\E ((Z_{n+1}^{(M)} - Z_n^{(M)} )^2 \| \F_n) \leq \Var 
   (Z_{n+1} - Z_n \| \F_n) \one_{\tau_M > n} \leq Q_n \one_{\tau_M > n}$$
and so $Z_n^{(M)}$ converges almost surely and in $L^2$ to 
a finite limit, $C_M$.  On the event $\{ \sum_n Q_n < \infty \}$, 
$\tau_M$ will be infinite for sufficiently large $M$, and the sequence 
$\{ Z_n \}$ will converge to $C_M + \sum_n \Delta_n$.
$\Cox$

\noindent{\sc Proof of Theorem}~\ref{th urn 2}: Let
$$Z_n = {X_n \over c Y_n} - \log Y_n \, .$$
We wish to apply Lemma~\ref{lem semi} to $\{ Z_n : n \geq 0 \}$,
so we must compute $\Delta_n$ and $Q_n$.  It will turn out that
$\Delta_n = O(1/Y_n^2)$ and $Q_n = O(n / Y_n^3)$, so we are going
to need a preliminary lower bound on the growth rate of $Y_n$
in order to see that these are almost surely summable.  

\begin{lem} \label{lem prelim}
For any $\ee > 0$, the function $X_n / Y_n^{1 + \ee}$ is a 
supermartingale when $X_n$ and $Y_n$ are both at least $c+2$.
It follows that $Y_n$ is almost surely eventually greater than
any power of $n$ less than 1.
\end{lem}

\noindent{\sc Proof:} To see that $X_n / Y_n^{1 + \ee}$ is a
supermartingale we compute the expected increment.
\begin{eqnarray*}
&& \E \left ( {X_{n+1} \over Y_{n+1}^{1 + \ee}} - {X_n \over Y_n^{1 + \ee}} 
   \| \F_n \right ) \\[2ex]
& = & {X_n \over X_n + Y_n} {1 \over Y_n^{1 + \ee}}
   + {Y_n \over X_n + Y_n} {c \over Y_n^{1 + \ee}} - {Y_n \over X_n + Y_n}
   {X_n ((Y_n + 1)^{1 + \ee} - Y_n^{1 + \ee}) \over Y_n^{1 + \ee} 
   (Y_n + 1)^{1 + \ee}} \\[2ex]
& = & {1 \over (X_n + Y_n) Y_n^{1 + \ee}} \left ( X_n + c Y_n - X_n 
   {Y_n - (Y_n / (Y_n + 1))^{1 + \ee}} \right ) \, .
\end{eqnarray*}
This is nonpositive when $\min \{ X_n , Y_n \} \geq c + 2$, proving
that $X_n / Y_n^{1 + \ee}$ is a supermartingale under this condition.
By Lemma~\ref{lem to infty}, both $X_n$ and $Y_n$ converge to infinity,
so there is an almost surely finite $N = N(\ee , \omega)$ such that
$\min \{ X_n , Y_n \} \geq c + 2$ for $n \geq N$, and consequently,  
$\{ X_n / Y_n^{1 + \ee} : n \geq N \}$ is a supermartingale.  This
implies that $\limsup_n X_n / Y_n^{1 + \ee} < \infty$, and hence 
for any $0 < \ee < \delta$, that $\limsup_n X_n^{(1 + \delta)^{-1}} / 
Y_n = 0$, proving the lemma.   $\Cox$

We continue with the proof of Theorem~\ref{th urn 2}.   
We first compute the expected increment $\Delta_n := \E (Z_{n+1} - Z_n 
\| \F_n )$ of $Z_n$.  
\begin{eqnarray}
\Delta_n & = & {X_n \over X_n + Y_n} {1 \over c Y_n}
   + {Y_n \over X_n + Y_n} {c \over c Y_n} + {Y_n \over X_n + Y_n}
   \left ( - {X_n \over c Y_n (Y_n + 1)} - \log {Y_n + 1 \over Y_n} 
   \right ) \nonumber \\[2ex]
& = & {1 \over X_n + Y_n} \left ( {X_n \over c Y_n } + 1 - {X_n \over
   c (Y_n + 1)} - Y_n \log (1 + 1/Y_n) \right ) \nonumber \\[2ex]
& = & {1 \over X_n + Y_n} \left ( {X_n \over c Y_n (Y_n + 1)} 
   + O ( {1 \over Y_n} ) \right ) \nonumber \\[2ex]
& = & O ( { 1 \over Y_n^2} ) . \label{eq delta}
\end{eqnarray}
Now compute an upper bound for the quadratic variation 
$Q_n := \E ((Z_{n+1} - Z_n)^2 \| \F_n)$ as follows.
\begin{eqnarray}
Q_n & = & {X_n \over X_n + Y_n} {1 \over c^2 Y_n^2} + {Y_n \over X_n + Y_n}
   \left ( {1 \over Y_n} - {X_n \over c Y_n (Y_n + 1)} - \log (1 +
   {1 \over Y_n}) \right )^2 \nonumber \\[2ex]
& \leq & {1 \over c^2 Y_n^2} + {1 \over Y^2} + \log^2 (1 + {1 \over Y_n}) 
   + {Y_n \over X_n + Y_n} {X_n^2 \over c^2 Y_n^4} \nonumber \\[2ex]
&\leq & {n \over C Y_n^3} \label{eq Q}
\end{eqnarray}
for an appropriate constant $C$, using the fact that $X_n 
\leq (1+c) n + X_0$.  We are now done: Lemma~\ref{lem prelim} together
with~(\ref{eq delta}) and~(\ref{eq Q}) show that $\Delta_n$ and $Q_n$
are almost surely summable, hence the conclusion of the theorem
follows from Lemma~\ref{lem semi}.   $\Cox$

\noindent{\sc Proof of Theorem~\ref{th urn 1}}:
We now prove Lemmas~\ref{lem supmart} and~\ref{lem general}, which 
together imply as a special case a result in the spirit of 
Theorem~\ref{th urn 1}.  We use supermartingales similar
to those in Lemma~\ref{lem prelim}.  

\begin{lem} \label{lem supmart}
Let $(X_n , Y_n)$ be a positive process converging coordinatewise to 
infinity.  Fix any $\beta > 1$ and suppose there is an $M = M(\beta) 
\leq \infty$ such that $Y_n^\beta / X_n$ is a supermartingale once 
$X_n , Y_n \geq M$.  Then 
$$\limsup_n {\log Y_n \over \log X_n} \leq {1 \over \beta} 
   \mbox{ on } \{ M < \infty \} \, .$$
Similarly, if $X_n / Y_n^\beta$ is a supermartingale once
$X_n , Y_n \geq M'$, then
$$\liminf_n {\log Y_n \over \log X_n} \geq {1 \over \beta} 
   \mbox{ on } \{ M < \infty \} \, .$$
\end{lem}

\noindent{\sc Proof}: Let $\tau_m$ be the least $n \geq m$ 
for which $\min \{ X_n , Y_n \} < M$.  Then 
$$\{ Y_{n \wedge \tau_m}^\beta /  X_{n \wedge \tau_m} : n \geq m \}$$
is a nonnegative supermartingale, so converges almost surely to a
limit $L(m)$. When $\tau_m = \infty$, it follows that $L(m)$ is the almost
sure limit of $Y_n^\beta / X_n$.  Thus when $\tau_m = \infty$,
$Y_n = [((L(m) + o(1)) X_n]^{1 / \beta}$.  If $L(m) > 0$ this implies 
$\limsup \log Y_n / \log X_n = 1 / \beta$, while if $L(m) = 0$,
the lim sup may be strictly less than $1 / \beta$.  On the event
$\{ M < \infty \}$ an $m$ exists with $\tau_m = \infty$, which finishes 
the proof of the first assertion.  The proof of the second assertion 
is similar.     $\Cox$

\begin{lem} \label{lem general}
Let $(X_n , Y_n)$ be a process generalizing the urn process in
Theorem~\ref{th urn 1} as follows.  For each $n$, with probability
$X_n / (X_n + Y_n)$, there is a $W > 0$ such that $X_{n+1} = W + X_n$ 
and $Y_{n+1} = Y_n$; with probability $Y_n / (X_n + Y_n)$, we have
$X_{n+1} = X_n$ and $Y_{n+1} = Y_n + 1$.  Suppose further that
with probability 1,
\begin{equation} \label{eq W1}
\E (W \| \F_n , X_{n+1} > X_n) \in [a,b]
\end{equation}
and
\begin{equation} \label{eq W2}
\E (W^2 \| \F_n , X_{n+1} > X_n) \leq K
\end{equation}
for some positive constant $K$ and some $0 < a \leq b$.  Then for
any $0 < \beta < a$, there is an $M$ such that whenever
$Y_n \geq M$, the function $Y_n^\beta / X_n$ is a supermartingale.
Likewise, for any $\beta > b$ there is an $M'$ such that
$X_n / Y_n^\beta$ is a supermartingale whenever $Y_n \geq M'$.
\end{lem}

\noindent{\sc Proof:} By a Taylor expansion, there exist constants
$c_1$ and $c_2$ such that for any $w$, and any sufficiently large $x$, 
$${1 \over x + w} \leq {1 \over x} - {w \over x^2} + c_1 {w^2 \over x^3 
   \, .}$$ 
Also, 
%
%
$$(y + 1)^\beta \leq y^\beta + \beta y^{\beta - 1} + c_2 \beta y^{\beta - 2} .$$
%
%
The expected increment $\Delta_n := Y_{n+1}^\beta / X_{n+1} 
- Y_n^\beta / X_n$, conditional on $\F_n$, is given by
$${X_n \over X_n + Y_n} \E \left ( {Y_n^\beta \over X_n + W} - 
   {Y_n^\beta \over X_n} \right ) + 
  {Y_n \over X_n + Y_n} \E \left ( {(Y_n + 1)^\beta \over X_n} - 
   {Y_n^\beta \over X_n} \right ) \, .$$
Plugging in the Taylor estimates above yields
$$\E \Delta_n \leq {X_n \over X_n + Y_n} {Y_n^\beta \over X_n^2} 
   \left ( - \E W + c_1 {\E W^2 \over X_n} \right ) +
  {Y_n \over X_n + Y_n} {\beta Y_n^{\beta - 1} \over X_n} 
   \left (1 + {c_2 \over Y_n} \right ) \, .$$
The assumptions on $W$ imply that
$$\E \Delta_n \leq {Y_n^\beta \over (X_n + Y_n) X_n} \left ( \beta
   - a + {c_2 \over Y_n} + {c_1 K \over X_n} \right ) \, .$$
When $\beta < a$ and $X_n$ and $Y_n$ are sufficiently large,
then this is nonpositive.  Choosing $M$ large enough so that
the constants $c_1$ and $c_2$ in the Taylor expansion are
valid whenever $X_n , Y_n \geq M$, we have proved the
first assertion of the lemma.

The proof of the second assertion is similar.  Choose $c$ so that
$$(y+1)^{- \beta} \leq y^{- \beta} - \beta y^{- \beta - 1} 
   (1 - {c \over y})$$
whenever $y \geq 1$.  
%
The expected increment $\Delta_n := 
X_{n+1} / Y_{n+1}^\beta  - X_n / Y_n^\beta$, conditional on $\F_n$, 
is given by
$${X_n \over X_n + Y_n} \E {W \over Y_n^\beta} +
%
%
  {Y_n \over X_n + Y_n} \E \left ( {X_{n} \over (Y_n + 1)^\beta} - 
   {X_n \over Y_n^\beta} \right ) \, .$$
Thus
$$\E \Delta_n \leq {X_n \over (X_n + Y_n) Y_n^\beta} \left (
   b - \beta + {\beta c \over Y_n} \right ) \, ,$$
proving the lemma for $M' = \beta c / (b - \beta)$.   $\Cox$

Finally, we show Lemma~\ref{lem supmart} and Lemma~\ref{lem general}
together imply the first order of approximation in Theorem~\ref{th urn 1},
namely that $\log X_n / \log Y_n \rightarrow a$ almost surely.
The urn process in Theorem~\ref{th urn 1}
satisfies the conditions of Lemma~\ref{lem general} with $K = a^2$ 
and $[a,b] = \{ a \}$.  Thus for any $0 < \beta < a$, we may plug
the conclusion of Lemma~\ref{lem general} into Lemma~\ref{lem supmart} 
to see that
$$\limsup_n {\log Y_n \over \log X_n} \leq {1 \over \beta} \, .$$
Similarly, for any $\beta > a$, we plug the conclusion of 
Lemma~\ref{lem general} into Lemma~\ref{lem supmart} to see that
$$\liminf_n {\log Y_n \over \log X_n} \geq {1 \over \beta} \, .$$
Since $\beta$ may be chosen arbitrarily close to $a$, we see that
$\log Y_n / \log X_n \rightarrow 1 / a$ as $n \rightarrow \infty$.
$\Cox$

\setcounter{equation}{0}

\section{Proof of Theorem~{\protect{\ref{th 3}}}}

In the next section we will prove the following lemma:
\begin{lem} \label{lem dominates}
There is an $\ee > 0$ such that for all integers $m > 0$,
$$\P (m+3 \in R \| m \in R) \leq 1 - \ee .$$
\end{lem} 
The first statement of Theorem~\ref{th 1} follows from Theorem~\ref{th 3}.
The second statement follows directly from Lemma~\ref{lem dominates}:
by induction, $\P (3n \in R) \leq (1 - \ee)^n$, which goes to
zero as $n \rightarrow \infty$.  Hence we concentrate on the proof 
of Theorem~\ref{th 3}.  Begin with a lemma.

\begin{lem} \label{lem tfae}
Let $J = [a , b]$ be an interval of integers containing zero.  Let $\P$ 
denote the law of VRRW on $\Z$ as before, and let $\P_J$ denote the law 
of a VRRW on the interval $J$, both started from 0.  Then the following
four conditions are equivalent.
\begin{quote}
$(i)$ $\P (R \subseteq J) > 0$; \\[1ex]
$(ii)$ $\P (R' \subseteq J) > 0$; \\[1ex]
$(iii)$ $\P \left ( \sum_n \one_{X_n = a} Z(n , a+1)^{-1} + 
   \sum_n \one_{X_n = b} Z(n , b-1)^{-1} < \infty \right ) > 0$. \\[1ex]
$(iv)$ $\P_J \left ( \sum_n \one_{X_n = a} Z(n , a+1)^{-1} + 
   \sum_n \one_{X_n = b} Z(n , b-1)^{-1} < \infty \right ) > 0$.
\end{quote}
\end{lem}

\noindent{\sc Proof:} We define a coupling, i.e., a measure $Q$ on 
pairs of paths $(\{ X_n : n \geq 0 \} , \{ X_n' : n \geq 0 \})$
such that the first coordinate of $Q$ has law $\P$ and the second 
has law $\P_J$.  To do so, choose $\{ X_n \}$ according to $\P$ and
let $\tau$ be the first time $n$ that $X_n \in \{ a-1 , b+1 \}$.
Let $X_n' = X_n$ for $n < \tau$, let $X_\tau' = a+1$ if $X_\tau = a-1$,
let $X_\tau' = b-1$ if $X_\tau = b+1$, and let $X_n'$ be chosen
from the transition probabilities for $\P_J$ independently of
$\{ X_n : n \geq 0 \}$ when $n > \tau$.  

Observe that 
$$Q(\tau = n+1 \| \tau > n) = \one_{X_n' = a} Z(n , a+1)^{-1} + 
   \one_{X_n' = a} Z(n , a+1)^{-1} .$$
Thus by Borel-Cantelli, $Q(\tau = \infty) > 0$ if and only
if condition~$(iv)$ is satisfied.  The event $\{ \tau = \infty \}$ 
is the same as the event $\{ R \subseteq J \}$, proving the 
equivalence of~$(i)$ and~$(iv)$.  Similarly, from the equation
$$Q(\tau = n+1 \| \tau > n) = \one_{X_n = a} Z(n , a+1)^{-1} + 
   \one_{X_n = a} Z(n , a+1)^{-1}$$
one sees that~$(i)$ and~$(iii)$ are equivalent.  The implication 
$(i) \Rightarrow (ii)$ is clear.  Finally, to see that~$(ii)$ 
implies~$(iii)$, assume~$(ii)$.  Thus with positive probability,
$Z(n , a-1) + Z(n , b+1)$ is bounded as $n \rightarrow \infty$.
By Borel-Cantelli, this means that 
$$ \P \left [ \sum_n \P (X_n \in \{ a-1 , b+1 \} \| \F_n ) < 
   \infty \right ] > 0 \, .$$
This sum is an upper bound for the sum in~$(iii)$, hence the
sum in~$(iii)$ is finite with positive probability.   $\Cox$

\begin{cor} \label{cor parts}
Suppose that $Z(n , a+1)$ and $Z(n , b-1)$ are $\Theta (n)$,
i.e., $\liminf Z(n , a+1) / n > 0$ and $\liminf Z(n , b-1) / n > 0$.
Then $\P (R \subseteq J) > 0$ if and only if 
$$ \P_J \left ( \sum_n {Z(n,a) + Z(n,b) \over n^2} < \infty \right ) > 0 .$$
\end{cor} 

\noindent{\sc Proof:} Let $\sigma_m$ be the first $n$ for which
$Z(n,a+1) = m$ and let $\rho_m$ be the first $n$ for which
$Z(n,b-1) = m$.  Summing by parts gives
\begin{eqnarray*}
&&\sum_k \one_{X_k = a} Z(k , a+1)^{-1} + \sum_k \one_{X_k = b} 
   Z(k , b-1)^{-1} \\[2ex]
& = & \sum_n {Z (\sigma_{n+1} , a) - Z(\sigma_n , a) \over n} +
   {Z (\rho_{n+1} , b) - Z(\rho_n , b) \over n} \\[2ex]
& = & \sum_n {Z(\sigma_n , a) - 1 \over n^2 - n} + {Z(\rho_n , b) - 1 
   \over n^2 - n} \, .
\end{eqnarray*}
Since $Z(n,r)$ is increasing in $n$ for all $r$ and we have assumed
$\sigma_n , \rho_n = O(n)$, this proves the corollary.   $\Cox$

\noindent{\sc Proof of Theorem}~\ref{th 3}: There are four steps to
the proof.  The first is to reduce to a VRRW on the five points
$-2, -1, 0, 1$ and 2.  The second is to show that this VRRW can, 
with positive probability, have $2 Z(n,1) / n$ remain in the 
interval $I$, while simultaneously $Z(n,2)$ and $Z(n,-2)$ remain
less than $n^{1 - \ee}$ for a prescribed $\ee = \ee (I) > 0$.  
The third step is to show that when these two things happen, then
actually $2 Z(n,1) / n$ converges to some $\alpha \in I$.  The fourth
step is to see that whenever $2 Z(n,1) / n$ converges, then
$Z(n,2)$ almost surely obeys the power law 
$$ \lim_{n \rightarrow \infty} {\log Z(n,2) \over \log n} =
   \lim_{n \rightarrow \infty} {2 Z(n,1) \over n} \, .$$

\ul{Step 1.}  This step is essentially done.  If we show that for
$J = \{ k-2 , k-1 , k , k+1 , k+2 \}$, the $\P_J$ probability of
properties $(ii)$~-~$(vi)$ holding simultaneously is positive,
then the conclusion of the theorem follows from $(ii)$~-~$(v)$
and Corollary~\ref{cor parts}.  The argument is the same for every
$k$, so from now on we assume without loss of generality that $k=0$, 
and set about proving Theorem~\ref{th 3} for $\P_J$ in place of $\P$,
where $J = \{ -2 , -1 , 0 , 1 , 2 \}$.  

\ul{Step 2.}
For the remainder of the argument, fix an interval $I = [c,d] \subseteq
(0,1)$ and a positive $\ee \leq \min \{ c , 1-d , d-c \} / 10$.  Let
$\beta = (1 - \ee)^{-1}$.
Also fix an integer $N_0$ and define stopping times depending on 
$N_0$ as follows.  Let $\tau_1$ be the least $n \geq N_0$ such
that $2 Z(n , 1) / n \notin I$.  Let $\tau_2$ be the least $n \geq N_0$
such that $Z(n , 2) \geq n^{1 - \ee}$ and let $\tau_3$ be the least
$n \geq N_0$ such that $Z(n , -2) \geq n^{1 - \ee}$.  Let 
$\tau = \tau_1 \wedge \tau_2 \wedge \tau_3$.  
Let $\{z_i : -2 \leq i \leq 2 \}$ be a quintuple of integers.  Our 
goal in this step is to identify an $N_0$ and a quintuple $z_i$ such that 
\begin{equation} \label{eq step 2}
\P_J (\tau = \infty \| Z(N_0 , i) = z_i : -2 \leq i \leq 2) > 0 \, ;
\end{equation}
in fact we will show it is near 1.  We assume~(\ref{eq step 2}) 
for the moment, and continue with steps~3 and~4.

\ul{Step 3.}  Let $\kappa_n$ be the time of the $n^{th}$ return
to the state 0 and define 
$$V_n = {Z(\kappa_n , 1) \over Z(\kappa_n , 1) + Z(\kappa_n , -1)}\; .$$
We will see below that when $\kappa_n < \tau_2
\wedge \tau_3$,
\begin{eqnarray} 
|\E_J (V_{n+1} - V_n \| \F_{\kappa_n}) | & \leq & C \kappa_n^{-1-\ee} 
   \label{eq step 3a} ; \\
\E_J ((V_{n+1} - V_n)^2 \| \F_{\kappa_n}) & \leq & C' \kappa_n^{-2} 
   \label{eq step 3b} .
\end{eqnarray}
Plugging these two bounds into Lemma~\ref{lem semi} shows that 
whenever $\tau = \infty$, the sequence $V_n$ must converge, 
to a value necessarily in $I$. This gives us parts~$(iv)$ and~$(v)$ of
the theorem, with part~$(vi)$ already following from step 2.

\ul{Step 4.}  We claim that for fixed $r$ and $s$, whenever 
\begin{equation} \label{eq st}
r \leq \liminf V_n \leq \limsup V_n \leq s
\end{equation}
and $\tau = \infty$, then
$${1 \over s} \leq \liminf {\log Z(n,2) \over \log n} \leq \limsup
   {\log Z(n , 2) \over \log n} \leq {1 \over r} \, .$$
To prove the claim, define the return times $\{ \alpha_n : n \geq 0 \}$ 
to state 1 by letting $\alpha_n = \min \{ n > \alpha_{n-1} : X_n = 1 \}$, 
and $\alpha_{-1}$ is set equal to $N_0 - 1$, for some $N_0$  Let 
\begin{eqnarray*}
U_n & = & Z(\alpha_n , 0) ; \\
U_n' & = & Z(\alpha_n , 2)  ;\\
\end{eqnarray*}
For any $\delta > 0$, we show that the conditions of 
Lemma~\ref{lem general} are satisfied with $(X_n , Y_n) = (U_n , U_n')$ 
and $[a,b] = [r - \delta , s + \delta]$.  Indeed, between times 
$\alpha_n$ and $\alpha_{n+1}$,
VRRW will either visit state 2 once or will visit state 0 some number
of times $W \geq 1$.  The probabilities of these disjoint cases are
respectively $U_n / (U_n + U_n')$ and $U_n' / (U_n + U_n')$.  Let $N_1$
be the least $N \geq N_0$ such that $r - \delta \leq \inf_{n \geq N}
V_n \leq \sup_{n \geq N} V_n \leq s + \delta$; when~(\ref{eq st}) holds,
$N_1$ will be finite.  We need to show that when $n \geq N_1$, 
then~(\ref{eq W1}) and~(\ref{eq W2}) hold, with $[a,b] = 
[r - \delta , s + \delta ]$.  Since $\P_J (W \geq k)$ is equal 
to the probability that on the first $k-1$ visits to state 0 
after time $\alpha_n$ the VRRW moves to the left, the assumption
that $n \geq N_1$ implies that 
$$ (1 - s - \delta)^k \leq \P_J (W \geq k) \leq (1 - r + \delta)^k ,$$
which gives  $1 / (s + \delta) \leq \E_J W \leq 1/(r - \delta)$ 
and $\E_J W^2 \leq K$ for some constant $K=K(r,\delta)$.  
The conclusion of Lemma~\ref{lem general} is
that $Y_n^{1/(s + 2 \delta)} / X_n$ and $X_n / Y_n^{1 / (r - 2 \delta)}$
are supermartingales for $n \geq M$, where $M$ will be finite
when~(\ref{eq st}) holds.  We then apply Lemma~\ref{lem supmart} 
together with the fact that $\alpha_n = O(n)$ on $\{ \tau
= \infty \}$ to see that
$${1 \over s + 2 \delta} \leq \liminf {\log Z(n , 2) \over \log n}
   \leq \limsup {\log Z(n , 2) \over \log n} \leq {1 \over r - 2 \delta}$$
on~(\ref{eq st}) when $\tau = \infty$.  Sending $\delta$ to 0 proves
the claim.  

Applying the claim simultaneously to all intervals $(r,s)$ with rational
endpoints, we see that conclusion~$(ii)$ of Theorem~\ref{th 3}
holds with probability 1 whenever $\tau = \infty$.  An identical
argument establishes conclusion~$(iii)$.  Since we have
shown that $\P_J (\tau = \infty)$ may be made arbitrarily close
to 1 by suitable choice of $\{ z_i : -2 \leq i \leq 2 \}$, we are done
with all four steps, modulo the verification of~(\ref{eq step 2}),
and of~(\ref{eq step 3a}) and~(\ref{eq step 3b}).

\ul{Cleanup step.}  First we prove~(\ref{eq step 3a}) 
and~(\ref{eq step 3b}).  Let $\Delta_n = V_{n+1} - V_n$.  We estimate 
$\Delta_n$ in three pieces.  Let $A$ be twice the least integer greater 
than $2 / \ee$.  Write $\Delta_n = R_n + S_n + T_n$ where
\begin{eqnarray*}
R_n & = & {Z (\kappa_n + 2 , 1) \over Z(\kappa_n + 2 , 1) + 
   Z(\kappa_n + 2 , -1)} - V_n \, , \\[1ex]
S_n & = & {Z (\kappa_{n+1} \wedge (\kappa_n + A) , 1) \over 
   Z (\kappa_{n+1} \wedge (\kappa_n + A) , 1) +
   Z (\kappa_{n+1} \wedge (\kappa_n + A) , -1) } - V_n - R_n \, , \\[1ex]
T_n & = & \Delta_{n} - R_n - S_n .
\end{eqnarray*}
By Proposition~\ref{pr polya}, $\E_J (R_n \| \F_{\kappa_n}) = 0$ 
and $R_n^2 \leq \kappa_n^{-2}$.  By the same token,
$S_n^2 \leq A^2 \kappa_n^{-2}$ and we easily see that
$$\E_J (|S_n| \| \F_{\kappa_n}) \leq {A \over \kappa_n} \P_J (\kappa_{n+1}
   > \kappa_n + 2 \| \F_{\kappa_n}) \leq {2 N_0^{1-\ee} 
   \over N_0 - 2 N_0^{1-\ee} }
   {A \over \kappa_n^{1 + \ee}}$$
when $\kappa_n < \tau_2 \wedge \tau_3$  
(as $n/2-n^{1-\ee}\leq\kappa_n\leq n$ and $Z(\kappa_n,0)= \kappa_n$).
Finally, since $T_n \leq 1$, 
we have
$$\E_J (|T_n|^i \| \F_{\kappa_n}) \leq \P_J (\kappa_{n+1} > \kappa_n
   + A \| \F_{\kappa_n})$$
for $i = 1, 2$.  The RHS is just the probability of at least
$A/2$ successive moves from state 1 to state 2 or state $-1$ to
state $-2$.  This probability is at most the maximum of
$$\prod_{i \leq A/2} {Z(\kappa_n , 2) + i \over Z(\kappa_n , 2)
   + Z(\kappa_n , 0) + i}$$
and the same expression with 2 replaced by $-2$.  Since 
$(\kappa_n^{-\ee})^{A/2} < \kappa_n^{-2}$ by choice of $A$,
and since the terms in the product are at most a constant multiple
of $\kappa_n^{-\ee}$ by the assumption that $\kappa_n < \tau_2 \wedge
\tau_3$, the RHS is bounded by a constant multiple of $\kappa_n^{-2}$.
Having bounded the conditional expectations of $R_n^2 , S_n^2$ and $T_n^2$ 
by multiples of $\kappa_n^{-2}$ and the magnitudes of
the conditional expectations of $R_n , S_n$ and $T_n$ by 
constant multiples of $\kappa_n^{-1-\ee}$, we have 
established~(\ref{eq step 3a}) and~(\ref{eq step 3b}).  

To establish~(\ref{eq step 2}), we will show that all of the
three probabilities $\P_J (\tau_2 \leq \tau_1 < \infty
\| \F_{N_0}))$, $\P_J (\tau_3 \leq \tau_1 < \infty \|
\F_{N_0})$ and $\P_J (\tau_1 < \tau_2 \wedge \tau_3 \| \F_{N_0})$ are 
simultaneously small when the values $Z({N_0},i) = z_i$ are chosen
appropriately.  As in step~4, we define the return times
$\alpha_n$ to state 1 by $\alpha_{-1} = N_0 - 1$ and $\alpha_{n+1} 
= \min \{ k > \alpha_n : X_k = 1 \}$.  Again set
$U_n = Z(\alpha_n , 0)$ and $U_n' = Z(\alpha_n , 2)$.  As in step~4,
the process $(U_n , U_n')$ evolves as the urns in Theorem~\ref{th urn 1},
where again we let $W = U_{n+1} - U_n$.  Assume that 
$\alpha_n < \tau_1 \wedge \tau_2$.  A lower bound for the
probability that $W > K$ is the probability that from state~1 the
VRRW visits 0 and then visits states~$-1$ and~0 $K$ times in alternation.   
Thus
\begin{eqnarray*}
\E_J (W \| \F_{\alpha_n} , W > 0) & \geq & 1 + \sum_{i = 2}^K
   \P_J (W > i - 1 \| \F_{\alpha_n} , W > 0) \\[1ex]
& \geq & 1 + \sum_{i=2}^K \prod_{j=1}^{i-1} {Z(\alpha_n , -1) + j - 1 
   \over Z(\alpha_n , -1) + Z(\alpha_n , 1) + j - 1} \; .
\end{eqnarray*}
There is a $\phi (K)$ such that when $N_0 \geq \phi (K)$ and
$\alpha_n < \tau_1 \wedge \tau_2$ then each factor in the product 
is at least $1 - d - 2 \ee$.  Thus for sufficiently large $K$ and
$\alpha_n < \tau_1 \wedge \tau_2$, the RHS is at least
$1 / (d + 3 \ee)$.  It is trivial to see that $\E_J (W^2 \| \F_n)$
is bounded.  Thus setting $a = 1 / (d + 3\ee)$ and
$\beta = 1 / (d + 4 \ee)$, we apply Lemma~\ref{lem general}
to see that $(U_n')^\beta / U_n$ is a supermartingale when    
$\alpha_n < \tau_1 \wedge \tau_2$.  More formally, let $\rho$
be the least $n$ for which $\alpha_n \geq \tau_1\wedge \tau_2$.
Setting $Y_n = U_{n \wedge \rho}'$ and $X_n = U_{n \wedge \rho}$, the
process $Y_n^\beta / X_n$ is a supermartingale.  Since $d + 4 \ee < 1 - \ee$,
we see that by definition of $\tau_1$ that $Y_n^\beta / X_n > 1$ 
if $\tau_2 \leq \tau_1 < \infty$, where $n$ is the least $j$ for which 
$\alpha_j \geq \tau_1$.  Therefore, by the supermartingale optional 
stopping theorem we arrive at
$$\P_J (\tau_2 \leq \tau_1 < \infty \| \F_{N_0}) \leq 
   {z_2^{1/(d+4\ee)} \over z_0} \; .$$

An entirely analogous argument with the states $-2$ and $-1$ in place
of 2 and 1 and $c$ in place of $1-d$ yields the analogous bound
$$\P_J (\tau_3 \leq \tau_1 < \infty \| \F_{N_0}) \leq 
   {z_{-2}^{1/(1-c-4\ee)} \over z_0 } \; .$$

Finally, we need to see how to make $\P_J (\tau_1 < \tau_2 \wedge \tau_3
\| \F_{N_0})$ small.  Let $\rho$ be the least $n$ for which $\kappa_n
\geq \tau_2 \wedge \tau_3$.  Then by~(\ref{eq step 3a}),
$$| \E W_\rho - W_0| \leq \sum_{m \geq \kappa_n} C \kappa_n^{-1-\ee}
   \leq C_1 (\ee) N_0^{- \ee} \, .$$ 
Similarly,~(\ref{eq step 3b}) gives
$$\Var (W_\rho - W_0) \leq \sum_{m \geq \kappa_n} C' \kappa_n^{-1}
   \leq C_1' N_0^{-1} \, .$$
On the other hand, if $\tau_1 < \tau_2 \wedge \tau_3$ then
$$|W_\rho - W_0| \geq \min \{ W_0 - c , d - W_0 \} \, .$$
Chebyshev's inequality applied to $W_\rho - W_0$ then shows that
$$\P_j (\tau_1 < \tau_2 \wedge \tau_3) \leq {C_1' N_0^{-1} \over
   (\min \{ W_0 - c , d - W_0 \} - C_1 N_0^{-\ee})^2 } \; .$$
When $N_0$ is sufficiently large, and $2 z_1 / N_0$ is sufficiently
close to $(c + d) / 2$, this is at most $C_2 N_0^{-1}$.  Thus
we have shown how to pick $z_{-2} , \ldots , z_2$ so that
$\P (\tau < \infty \| \F_{N_0})$ can be made arbitrarily small, 
which finishes the proof of~(\ref{eq step 2}) and of 
Theorem~\ref{th 3}.    $\Cox$

\section{Remaining proofs and open questions}

The proof of Lemma~\ref{lem dominates} is quite similar to the proof
of Theorem~\ref{th 3}.  We give an outline for the argument, leaving
out details that are the same as in the proof of Theorem~\ref{th 3}.

\noindent{\sc Sketch of proof of Lemma}~\ref{lem dominates}:  The first
reduction is to analyze VRRW on $(-\infty , m+2]$. 
Let $\tau_m$ be the first time $m$ is reached.  The hardest
part, because it requires a simultaneous induction on two different
stopping times, similar to~(\ref{eq step 2}), is the following:
\begin{quote} 
Claim 1: There is a constant $\delta > 0$ such that for all $m$,
the probability is at least $\delta$ that
inequalities $(i)$ and $(ii)$ hold for every $n \geq \tau_m$:
\begin{eqnarray*}
(i) && Z(n,m) \geq Z(n,m+2)^2 \, ; \\
(ii) && Z(n,m-1) \geq 2 Z(n,m+1) \, .
\end{eqnarray*}
\end{quote}
The other essential ingredient is 
\begin{quote}
Claim 2: $Z(n , m+1) \leq (1/4) Z(n , m+2)^2 \one_{A_n}$ finitely often
almost surely, where $A_M$ is the event that $(i)$ and $(ii)$
of the previous claim are true for all $n \in [\tau_m , M]$.  
\end{quote}
Assume these two claims and let $\sigma_k$ be the time of the $k^{th}$
visit to site $m+2$.  From the first claim, the decreasing limit
$A_\infty$ has probability at least $\delta$.  Since $Z(\sigma_k , m+2)
= k+1$, it follows from the second claim that on $A_\infty$,
$$\sum_{k=1}^\infty Z(\sigma_k , m+1)^{-1} < \infty . $$
Fix $M , \ee > 0$ such that with probability at least $\ee$,
$\sum_{k=1}^\infty Z(\sigma_k , m+1)^{-1} < M$.
As in Corollary~\ref{cor parts} it then follows for VRRW on $\Z$
that $\P (m+3 \notin R \| m \in R) > 0$ and in fact that a lower
bound is $\delta := \ee \exp (-2M)$.  

To prove the first claim, stop the walk the first time either
condition~$(i)$ or~$(ii)$ is violated.  Consider first the process
$\{ (U_n , V_n) \} := \{ (Z(\rho_n , m-1) , Z(\rho_n , m+1)) \}$, 
where $\rho_n$ are the successive hitting times of site $m$.  
At each step, precisely one of the coordinates is updated, with 
Polya-like probabilities, so by Proposition~\ref{pr polya}, 
the expected increment of $U_n / (U_n + V_n)$ is given by the
contributions from increments of magnitude greater than 1:
$$\E {U_{n+1} \over U_{n+1} + V_{n+1}} - {U_n \over U_n + V_n} = 
   \E {(U_{n+1} - U_n - 1)^+ - (V_{n+1} - V_n - 1)^+ \over
   U_{n+1} + V_{n+1} } \, .$$
The term $(U_{n+1} - U_n - 1)^+$ is nonnegative and the term
$$\E {-(V_{n+1} - V_n - 1)^+ \over U_{n+1} + V_{n+1}}$$ is of order 
$${1 \over U_n + V_n} {Z(\rho_n , m+2) \over Z(\rho_n , m)} = 
   O(n^{-1} n^{-1/2})$$
by condition~$(ii)$.  This expresses $U_n / (U_n + V_n)$ as a 
martingale plus a drift term whose negative part is summable.

Consider next the process
$\{ (U_n' , V_n') \} := \{ (Z(\rho_n , m) , Z(\rho_n , m+2)) \}$, 
where $\rho_n$ are now the successive hitting times of $m+1$.  Again
the updates are in a single coordinate chosen with P\'olya-like
probabilities, with the increment in $V_n'$ being 1 and the
increment in $U_n'$ having conditional mean at least 3.  Using
Lemma~\ref{lem general}, just as in Step~4 of the proof of 
Theorem~\ref{th 3}, we see that $(V_n')^2 / U_n'$ is a supermartingale.

The optional stopping theorem now shows that from an appropriate
initial position, the probability of stopping due to a violation
of~$(i)$ or~$(ii)$ is arbitrarily low.  The initial position (or
one at least as good) can be attained with a probability bounded
away from zero (unless $m$ is visited only finitely often, which
is even better!) so the claim is proved.  

Finally, to prove the second claim, let $X_n = Z(\rho_{n+1} , m+1)
- Z(\rho_n , m+1)$ where now $\rho_n$ are the successive hitting times
of site $m+2$.  On the event $A_{\rho_{n+1}}$, the probability of
a transition to $m+2$ from $m+1$ between times $\rho_n$ and 
$\rho_{n+1}$ is bounded above by $1/(n+1)$ (use condition~$(i)$ and
$Z(\cdot , m+2) = n$).  Thus $X_n$ stochastically dominates a 
geometric of mean $n$, and it is easy to verify that 
$$\sum_{k=1}^n X_k < {n^2 \over 4} \mbox{ finitely often, a.s.},$$
which proves the second claim and hence the lemma.    $\Cox$

\noindent{\sc Proof of Theorem~\ref{th 2}:} 
It is straightforward that $\P (|R'| = 2)=\P (|R'| = 3)  = 0$, so we
suppose that $R'=J:=\{-2,-1,0,1\}$ and show that this leads 
to a contradiction. 
To simplify notation, let $A_n$, $B_n$, $C_n$ and  $D_n$ denote
$Z(n,x)$ for $x=-2,-1,0$ and $1$ respectively. 
The assumption $R'=J$ yields that there exist $N$ such that
$X_n\in J$ as soon as  $n\geq N$. Throughout the rest of the proof, 
we assume that $n\geq N$.

Let $\kappa_m$ be the time of the $m^{th}$ return to $-1$
(clearly, $\kappa_m$'s are $<\infty$ on the event $R'=J$).
Between the times $\kappa_m$ and $\kappa_{m+1}$, the random walk
goes from $-1$ either to $-2$ and returns
to $-1$, or to $0$ and (possibly) bounces between $0$ and $1$ 
before going back to $-1$.
Therefore,
$W_m:=A_{\kappa_m}/(A_{\kappa_m}+C_{\kappa_m})$ 
is a non-negative supermartingale
which converges a.s. to some random variable $W$.
Consider two cases:
\begin{itemize}
\item $W(\omega)=w>0$;
\item $W(\omega)=0$.
\end{itemize}
In the first case, there exists $N_1>N_0$ 
such that $W_m>w/2$ for $n>N_1$.
As in Lemma~\ref{lem tfae}, the probability never to
jump from $-2$ to $-3$ is  positive only whenever
\begin{eqnarray*}
 \sum_n \frac{A_{n+1}-A_n}{B_n}<\infty
\end{eqnarray*}
which is equivalent to the following sum being finite
\begin{eqnarray}\label{eq nl4}
 \sum_{n\geq N_1} \frac{A_n}{B_n}\frac{B_n-B_{n-1}}{B_{n-1}}
 \geq \frac w2 \sum \frac{B_n-B_{n-1}}{B_{n-1}}
\end{eqnarray}
(here we used the obvious inequality $B_n\leq A_n+C_n$). 
However, the sum in the RHS of~(\ref{eq nl4})
is a tail of the harmonic series and therefore diverges.

Before we proceed to the second case we observe that by the same arguments
we can restrict the problem to the case when both
$A_n/(A_n+C_n)$ and $D_n/(D_n+B_n)$ go to zero.
As a result, $A_n/C_n\to 0$ and $D_n/B_n\to 0$ as well. 
Taking into account that $B_n\leq C_n+A_n$ and $C_n\leq B_n+D_n$ 
we conclude that
\begin{eqnarray}\label{eq nlas}
\frac{B_n}{C_n}\to 1,\ \frac{2B_n}{n}\to 1,\  \frac{2C_n}{n}\to 1. 
\end{eqnarray}

Let $\tau_m$ be the time of the $m^{th}$ visit
to $-1$ or $0$ skipping at least one step, 
$$
\tau_m:=\inf \{n>\tau_{m-1}+1:\ X_n\in\{-1,0\}\},
\ \ \ \tau_0:=N_0.
$$
Define
\begin{eqnarray*}
U_m=\frac{B_{\tau_m}+C_{\tau_m}}2,\ \ \ V_m=A_{\tau_m}+D_{\tau_m}.
\end{eqnarray*}
If $X_{\tau_m}=-1$ ($X_{\tau_m}=0$ resp.), then
between the times $\tau_m$ and $\tau_{m+1}$ VRRW will either
(1) go to the left  (right resp.) and back, or
(2) go to the right (left  resp.) and back, or
(3) go twice to the right (left  resp.) and make one step back. 
Consequently, either $U_{m+1}=U_m+1$ and $V_{m+1}=V_m$ (when (2) takes place)
or  $U_{m+1}\leq U_m+1$ and $V_{m+1}=V_m+1$ (when (1) or (3) takes place).
We claim that the probability of the latter event, denoted by $F$,
is greater than $V_m/(U_m+V_m)$ when $m$ is large enough.
This, in turn, will imply that the process $(U,V)$ can be coupled with 
some general urn model process described by (\ref{eq urn})
such that $U_m\leq X_m'$ and $V_m\geq Y_m'$.

To prove that $\P(F)\geq V_m/(U_m+V_m)$ we consider
the quantity $A_n-B_n+C_n-D_n$ which is ``almost'' invariant for $n\geq N$.
Namely, there exists a (possibly negative) constant $K$, 
depending on the history
of VRRW before time $N$ only, such that $A_n-B_n+C_n-D_n=K$ whenever $X_n=-1$
and $A_n-B_n+C_n-D_n=K+1$ whenever $X_n=0$. 
If we denote $t_m:=B_{\tau_m}+D_{\tau_m}$, then
$A_{\tau_m}+C_{\tau_m}$ equals $t_m+K$ or $t_m+K+1$ 
when VRRW is at $-1$ or at $0$ respectively. 
In the second case
$$
 \P(F)=\frac {D}{t}+\left(1-\frac{D}{t}\right)
  \frac{A}{t+K+1}= \frac{V-AD/(t+K+1)}{t} -\frac{A(K+1)}{t(t+K+1)}
$$
(we omit the indices for simplicty). Taking into account
that $2U+V=2t+K+1$, $V>A$ and $AD\leq V^2/4$,  we obtain
\begin{eqnarray*}
 \P(F)-\frac{V}{U+V}&\geq& \frac{V-V^2/(4t+4K+4)}{t} -\frac{A(K+1)}{t(t+K+1)}
 -\frac{V}{t+(V+K+1)/2}
\\&\geq&\frac{V(V+K-1)}{t(2t+V+K+1)}-
 \frac{V^2}{t(4t+4K+4)}-\frac{V(|K|+1)}{t(t+K+1)}
\end{eqnarray*}
As $m\to\infty$ (and, therefore, $n\to\infty$) we have $V_m\to\infty$
and $V_m=o(t_m)$, 
whence 
$$
\P(F)-\frac{V}{U+V}
=\frac{V^2}{4t^2}\left(1-\Theta\left(\frac Vt\right)
-\Theta\left(\frac 1V\right)\right)=\frac{V^2}{4t^2}\left(1-o(1)\right)
$$
is non-negative for $m\geq M$ where $M$ is some
constant. The case $X_{\tau_m}=-1$ can be analyzed
in the similar way.

We have shown that for large $m$ the process $(U_m,V_m)$ can be coupled
with the process $(X_m',Y_m')$ obeying the law
(\ref{eq urn}) with $a=d=1$, $b=0$ and $c=1$
such that  $U_m\leq X_m'$ and $V_m\geq Y_m'$.
By Theorem~\ref{th urn 2}, there exists
$$
 \lim_{m\to\infty} \frac{X_m'}{Y_m'} -\log(Y_m')\in(-\infty,\infty)
$$
which, in turn, implies the existence of a random variable 
 $\zeta\in(0,\infty)$ such that
$$
X_m'\leq Y_m'\log(2\zeta Y_m')\mbox{\ \ \  for all $m\geq M$.}
$$
Clearly, $X_m'/Y_m'\to \infty$, so
$$
Y_m'\geq \frac{X_m'}{\log(2\zeta X_m')}
$$
for all $m$ larger than some $M_1\geq M$.
Since $U_m\leq X_m'$, $V_m\geq Y_m'$ and the function $f(x)=x/\log(2\zeta x)$
is increasing for large $x$,  $V_m\geq U_m/\log(2\zeta U_m')$.
Furthermore, $\tau_m\to\infty$ and $\tau_{m+1}\leq \tau_{m}+3$, so  
we  asymptotically have
$$
A_n+D_n\geq \frac{(B_n+C_n)/2}{\log( \zeta (B_n+C_n) )}\simeq \frac{n}{2\log(n)}
$$
by (\ref{eq nlas}).

The event that VRRW does not jump off $J$ can occur
only when the sum
$$
\sum \frac {A_{n+1}-A_{n}}{B_n}+ \frac{D_{n+1}-D_{n}}{C_n}
$$
is finite. Summing by parts as in Corollary~\ref{cor parts}, 
we obtain that this is equivalent to the finiteness of the sum
$$
\sum_n \frac{A_n+D_n}{n^2}\geq
const+\sum_n \frac{(B_n+C_n)/2}{\log( \zeta (B_n+C_n) )}\frac{1}{n^2}
\simeq const+\sum_n \frac{1}{2n \log(\zeta n)}
$$
which diverges.
Therefore, $\P (R' = J) = 0$, completing the proof.   $\Cox$
%
%

We end with some questions.  The strongest conjecture about VRRW
on $\Z$ is the one stated after Theorem~\ref{th 3}, to the effect
that the behavior described in Theorem~\ref{th 3} happens with
probability 1.  Some smaller steps toward this would be to prove
that the set of sites visited with positive density must be 
connected and to prove that $\alpha$ can never be 0 or 1.  This 
would, for example, rule out that sites $-1$ and 0 are visited
with density $1/2$, and by time $n$ the numbers of visits to the sites 
$1 , 2 , \ldots$ are asymptotically $n / \log n , n / (\log n \log \log n) 
,\ldots$.  Another graph on which VRRW may have interesting behavior
is $\Z^2$.  
Ferrari and Meilijson (personal communication, 1996)
also have some results about VRRW on a tree. 

A further question is that of stochastically comparing VRRW with
different histories.  For example, we originally thought we could 
prove a version of Lemma~\ref{lem dominates} in which it was shown
that $\P (m+3 \in R \| m \in R) \leq \P (3 \in R)$, by showing
that the extra weight to the left of $m$ the first time $m$ is reached
can only help the range stay bounded above by $m+2$.  We were unable
to do this, by coupling or martingale arguments, but believe that
some such comparison must hold.  The easiest to state are false.

\newpage\noindent
Robin Pemantle\\
Department of Mathematics, University of Wisconsin-Madison\\
Van Vleck Hall, 480 Lincoln Drive, Madison, WI 53706\\

\medskip\noindent
Stanislav Volkov\\
The Fields Institute for Research in Mathematical Sciences\\
222 College Street, Toronto, Ontario, Canada, M5T3J1\\

\end{document}